\documentclass[11pt,reqno]{amsart}
\usepackage {amscd}

\newtheorem{theorem}{Theorem}[section]
\newtheorem{lemma}[theorem]{Lemma}
\newtheorem{prop}[theorem]{Proposition}
\newtheorem{kor}[theorem]{Corollary}
\newtheorem*{theoremho}{Main Theorem}

\newcommand{\sus}{\subseteq}
\newcommand{\od}{>}
\newcommand{\gin}[1]{\makebox{gin}(#1)}
\newcommand{\ini}[1]{\makebox{in}(#1)}
\newcommand{\ii}[1]{\makebox{in}(#1)}
\newcommand{\lpil}{\longrightarrow}

\newcommand{\ptre}{{\bf P}^3}
\newcommand{\bx}{{\mbox{\boldmath $x$}}}
\newcommand{\bt}{{\mbox{\boldmath $t$}}}
\newcommand{\go}{{\mathcal O}}

\newcommand{\pd}{\partial}

\newcommand{\saf}{{\bf A}^s}
\newcommand{\ba}[1]{\overline {#1}}

\newcommand{\ga}{\gamma}
\newcommand{\Ga}{\Gamma}
\newcommand{\si}{\circ}
\newcommand{\sas}{{\bf A}^{s-1}}
\newcommand{\pil}{\rightarrow}
\newcommand{\te}{\otimes}

\newcommand{\parg}{{\vskip 2mm \addtocounter{theorem}{1}  

                   \noindent {\bf \thetheorem} \hskip 1.5mm }}

\begin{document}

\title {A property deducible from the generic initial ideal}
\author { Gunnar Fl{\o}ystad}
\address{ Matematisk Institutt\\
          Allegaten 55 \\
          5007 Bergen }
\email{ gunnar@mi.uib.no}

\begin{abstract}
Let $S_d$ be the vector space of monomials of degree $d$ in the variables
$x_1, \ldots, x_s$. For a subspace $V \sus S_d$ which is in general 
coordinates, consider the subspace $\gin V \sus S_d$ generated by initial
monomials of polynomials in $V$ for the revlex order. We address the question
of what properties of $V$ may be deduced from $\gin V$.

This is an approach for understanding what algebraic or geometric 
properties of a homogeneous ideal $I \sus k[x_1, \ldots, x_s]$ that may
be deduced from its generic initial ideal $\gin I$.
\vskip 6mm
\noindent {\it 1991 Math. Subj. Class.:} 13P10

\end{abstract}

\maketitle

\section*{Introduction}

During the recent years the generic initial ideal of a homogeneous ideal
has attracted some attention as an invariant. 
An intriguing problem is what algebraic or geometric properties of the
original ideal can be deduced from the generic initial ideal.

In this paper we take perhaps the most elementary approach possible.
 Let 
$S = k[x_1, \ldots, x_s]$ and let $\od$ be the reverse lexicographic order
of the monomials in $S$. Denote by $S_d$ the graded piece of degree $d$ in 
$S$. Suppose $V \sus S_d$ is a subspace. Denote by $\gin V$ the subspace
of $S_d$ generated by initial monomials 
of polynomials of the subspace of $S_d$ obtained from $V$ by performing 
a general change of coordinates.
Then one may ask what properties of $V$ may be 
deduced from $\gin V$? The following result gives an insight in this
direction.

Let $W = (x_1, \ldots, x_r) \sus S_1$ which is a linear space. 
Suppose that $s \geq r \geq 3$. 

\begin{theoremho} Let $V \sus S_{n+m}$ be a linear space such that 
\[ \gin V = W^n x_1^m \sus S_{n+m}. \]
Then there exists a polynomial $p \in S_{m}$ and a linear subspace
$W_n \sus S_n$ such that $V = W_n p.$
\end{theoremho}

Note that if $s=r$ then $W^n x_1^m$ are the largest monomials in 
$S_{n+m}$ for the lexicographic order. Thus if $\od$ had been the
lexicographic order and $\gin V = W^n x_1^m$ then we could deduce 
virtually nothing about $V$.

\vskip 2mm

The general idea of the proof is inspired by Green and worth attention 
because of its seeming naturality in dealing with problems of this kind.

The idea in its vaguest and most generally applicable form is the following.
Suppose $\gin V$ has a given form, and suppose $V$ is in general coordinates
so $\ii V = \gin V$. 
The given form of $\ii V$ implies some algebraic or geometric property of
$V$. Let now $g : S_1 \pil S_1$ be a general change of coordinates.
Then $\ii {g^{-1}.V} = \gin V$ also. Thus $g^{-1}.V$ will also have this
property. Then this property may be translated back to a property of $V$.
This gives a continuous set of properties that $V$ will satisfy. From this one
may proceed making deductions about what $V$ may look like.

\vskip 2mm
In this paper this is applied concretely as follows. In the case $r=s$ the 
given form of $\ii V = \gin V$ implies that there is a $p_1 $ in $S_m$ such
that $ x_r^n \cdot p_1 \in V$.
The fact that $\ii {g^{-1}.V} = \gin V$ also implies that there is a 
$p_{g^{-1}}$ in $S_m$ such that $ x_r^n \cdot p_{g^{-1}} \in g^{-1}.V$.
Translating this property back to $V$ we get
\begin{equation} (g.x_r)^n \cdot g.p_{g^{-1}} \in V. 
\label{hnpz} \end{equation}
Now for the family of linear forms $h = \sum t_ix_i$ one may choose a general
family of $g$'s depending on $h$ such that $g.x_r = h$. Then 
equation (\ref{hnpz}) may be written as
\begin{equation} h^n p \in V \label{hnp0} \end{equation}
where $p$ is a form
of degree $m$ depending on $h$. 

\vskip 2mm The second technique, specifically suggested by Green,
is to {\it differentiate} this equation
with respect to the $t_i$. All the derivatives will still be in $V$.
(This is just the fact that when a vector varies in a vector space its
derivative is also in the vector space.)
Letting $V_{|h=0}$ be the image by the composition $V \pil S \pil S/(h)$ this 
enables us to show that the forms in $V_{|h=0}$ have a common factor
of degree $m$.

\vskip 2mm The third basic ingredient is now proposition 3.4 which says that
if the $V_{|h=0}$ have a common factor of degree $m$, then $V$ has a common
factor of degree $m$. 

Having proven the case $s=r$, the case $s > r$ 
may now be proven by an induction process.

\vskip 3mm

The organization of the paper is as follows. In the first three sections we 
develop general theory which does not presuppose {\it anything} about what
$\gin V$ actually is.

In section 1 we give some basic definitions and notions.
In section 2 we define the generical initial space of a subspace $V$ of $S$
by using a {\it generic} coordinate change on $V$. We also give some basic
theory for this setting which will be used in sections 4 and 5.

Section 3 presents the framework in which we will work. Instead of considering
a continuously varying form $h= \sum_{i=1}^s t_i x_i$ 
in $k[x_1, \ldots, x_s]$, we
consider $h$ as a linear form in $K[x_1, \ldots, x_s]$ where 
$K = k(t_1, \ldots, t_s)$, the field of rational functions of the $t_i$'s.

If now $V \sus k[x_1, \ldots, x_s]_d$ is a subspace let $V_K = V \te_k K
\sus K[x_1, \ldots, x_s]$. The main result here, proposition 3.4, says that if
the forms in $V_{K|h=0}$ have a common factor of degree $m$ 
then the forms in $V$ have a common factor
of degree $m$. This is proven using differentiation of forms with respect
to the $t_i$.

Only {\it from now on} do we assume that $\gin V$ has the special form
given in the main theorem.
In section 4 we prove the case $s=r$ in the main theorem.
Section 5 proves the case $s > r$ of the main theorem.
In section 6 we give an application of the main theorem. The example
originated in discussions with Green and was what triggered this paper.
Consider the complete intersection of three quadratic forms in $\ptre$.
Let $I \sus k[x_1,x_2,x_3,x_4]$ be its homogeneous ideal. By standard theory 
one may deduce that there are two candidates for $\gin I$ :

\begin{eqnarray*} 
J^{(1)} & = & (x_1^2, x_1x_2, x_2^2, x_1x_3^2, x_2x_3^2, x_3^4), \\
 J^{(2)} & = & (x_1^2, x_1x_2, x_1x_3, x_2^3, x_2^2x_3, 
x_2x_3^2, x_3^4).
\end{eqnarray*}

By the main theorem, if $\gin I = J^{(2)}$ then the quadratic forms
in $I_2 \sus S_2$ would have to have a common factor. Impossible.
Thus $\gin I = J^{(1)}$.

\vskip 2mm
Throughout the article all fields have characteristic zero.

\section{Basic definitions and notions}

\parg Let $S = k[x_1, \ldots, x_s]$.
The graded piece of degree $d$ is denoted by $S_d$.
If $I = (i_1, i_2, \ldots, i_s)$ we use the notation
\[ \bx^I = x_1^{i_1}\cdots x_s^{i_s}. \]
It has degree $|I| = \sum i_j$. 


Suppose now we have given a total order on the monomials.
For a homogeneous polynomial $ f = \sum a_I\bx^I$ in $S$
(henceforth often referred to as a form) let the initial monomial be
\[ \ii f = \max \{ \bx^I \, | \, a_I \neq 0 \}. \]

For a homogeneous vector subspace $V \sus S$ let the {\it initial subspace} be
\[ \ini V = ( \{ \ini f \, | \, f \in V\} ) \]
the homogeneous vector subspace of $S$ generated by the initial monomials
of forms in $V$. 

Sometimes we wish to consider another polynomial ring $R[x_1, \ldots, x_r]$
where $R$ is a commutative ring. Denote this by $S_R$.
  The initial monomials $\ii f$ for $f \in V$ may equally well be 
considered as elements of $S_R$. We may thus speak of $\ii V$ over $R$
(when $V \sus S$) which is the free $R$-module in $S_R$ generated by
$\{ \ini f \, | \, f \in V\}$.

\parg The monomial order we shall be concerned with in sections 4 and 5
is the reverse lexicographic order. 
Then the monomials of a given degree is ordered by
$\bx^I > \bx^J$ if  $i_r < j_r$ where $r$ is the greatest number
with $i_r \neq j_r$. Intuitively $\bx^J$ is "dragged down" by having a large
"weight in the rear".

\parg
For a linear form $l \in S_1$ denote by
$V|_{l=0}$ the image of the composition
\[ V \lpil S \lpil S/(l). \]
The following basic fact for the revlex order, 
proposition 15.12 a. in \cite{Ei}, will be used several times
\[ \ii {V_{|x_s = 0}} = \ii V _{|x_s = 0}. \]

 \section{The generic initial space}

   The following section contains the definition of the generic initial space
and some general theory related to it. The things presented here are 
certainly in the background knowledge of people but due to a lack of suitable
references for a proper algebraic treatment we develop the theory here.
The most important things are proposition 2.9 and paragraph 2.11.

 \parg We identify $S = k[x_1, \ldots, x_s]$ as the affine coordinate ring of 
$\saf$. Let $G = GL(S_1^{\vee})$. There is a natural action
\[ \saf \times G \lpil \saf \]
given by $(a,g) \mapsto g^{-1}.a$. This gives a $k$-algebra homomorphism
\[ \gamma : k[x_1, \ldots, x_s] \lpil k[x_1, \ldots, x_s] \te_k k[G]. \]
If $R$ is a $k[G]$-algebra, we also by composition obtain a $k$-algebra
homomorphism
\[ \ga_R : k[x_1, \ldots, x_s] \lpil k[x_1, \ldots, x_s] \te_k k[G]
    \lpil  R[x_1, \ldots, x_s]. \]
Note that if $R = k(g)$ for a point $g \in G$, then $\ga_{k(g)}$ is just
the action of $g$ on $k[x_1, \ldots, x_s]$.

Let $K_G$ be the function field of $G$. The image of a homogeneous 
subspace $V \sus S$ by $\ga_{K_G}$ generates a homogeneous subspace
$(\ga_{K_G}(V))$ of the same dimension as $V$. 
Suppose now a total monomial order is given.
The initial monomials of 
$(\ga_{K_G}(V))$ generate a linear subspace over $k$ (or over $K_G$), which is 
called the {\it generic initial subspace} of $V$ over $k$ (or over $K_G$)
and is denoted $\gin V$. Henceforth we shall drop the outer 
paranthesis of $(\ga_{K_G}(V))$ and write this as $\ga_{K_G}(V)$.

\parg Let $\gin V = (m_1, \ldots, m_t)$ for some monomials $m_i$.
Let $b_i \in \ga_{K_G}(V)$ be such that
\[ b_i = m_i + b_{i0} \]
where $b_{i0}$ consists of monomials less than $m_i$ for the given order.
Now there is an open subset $U \sus G$ such that all the $b_i$ lift to 
elements of $\go(U) [x_1, \ldots, x_s]$. Now we immediately get.

\begin{prop}

There is an open subset $U \sus G$ (take the one above) such that for $g \in U$
then
\[ \ii {(\ga_{k(g)}(V))} = \gin V \, (\mbox{over } k(g)).\]
\end{prop}

(The original reference for this is \cite{Ga}.)

\parg Now choose a $g_0 \in G$ such that $k(g_0) = k$. There is then a diagram
\[ \begin{CD}
 \saf_{K_G} @>\alpha_{g_0}>> \saf_{K_G} \\
 @VVV @VVV \\
 \saf \times G  @>{\cdot g_0}>> \saf \times G \\
 @VVV @VVV \\
 \saf @>{g_0^{-1}.}>> \saf.
\end{CD} \]
The lower horizontal map is the natural action. The middle map is given by
$(a,g) \mapsto (a, g g_0)$ and the lower vertical maps are just the action
of $G$. The upper horizontal map is the map induced by the middle map.
From the commutativity of the diagram we see that
\[ \ga_{K_G}(g_0.V) = \alpha_{g_0}^*( \ga_{K_G}(V)) \]
where $\alpha_{g_0}^*$ is the automorphism of $K_G[x_1, \ldots, x_s]$
induced by $\alpha_{g_0}$. Note that $\alpha_{g_0}^*$ comes from an
automorphism of $K_G$. 
So it does not affect the 
variables $x_i$. 

Thus we see that the $\alpha_{g_0}^*(b_i) = m_i + \alpha_{g_0}^*(b_{i0})$
are a basis for $\ga_{K_G}(g_0.V)$, where the monomials in 
$\alpha_{g_0}^*(b_{i0})$ are less then $m_i$ for the given order.
Also note that the $\alpha_{g_0}^*(b_i)$ lift to the open subset 
$U.g_0^{-1} \sus G$.
Thus we have proven the following.

\begin{lemma}

Given $g \in G$, by replacing the subspace $V$ by $g_0.V$ and the 
open subset $U$ by $U.g_0^{-1}$ for a suitable $g_0$, we may
assume that $g$ is in the open subset from proposition 2.3.

\end{lemma}

\parg Now let $\phi : X \lpil G$ be a morphism. We get a morphism
\[ \saf \times X \lpil \saf \]
and thus a $k$-algebra morphism
\[ \ga_{K_X} : k[x_1, \ldots, x_s] \lpil K_X[x_1, \ldots, x_s] \]
where $K_X$ is the function field of $X$. 
We get a homogeneous subspace $( \ga_{K_X}(V))$ and also here we shall
henceforth drop the outer paranthesis.
By performing a suitable
coordinate change of $V$ we may assume (by lemma 2.5) that 
$\phi(X) \cap U \neq \emptyset$.
The following is now immediate from the results above.

\begin{lemma}

\begin{itemize}
\item [1.] For $x$ in the open subset $\phi^{-1}(U) \sus X$ we have
\[ \ii {\ga_{k(x)}(V)} = \gin V \, (\mbox{over } k(x)). \]

\item [2.] $\ii {\ga_{K_X} (V)} = \gin V \, (\mbox{over } K_X).$ 

\item [3.] Given $x \in X$ then we may assume that $\phi(x) \in U$.

\end{itemize}

\end{lemma} 

\parg By 1.3 we have
\[ \ii {V_{|x_s = 0}} = \ii V _{|x_s = 0}. \]
We would like to have a suitable version of this for generic subspaces.
The version we need is 2. in the following. It is used most importantly
in the proof of lemma 5.2

\begin{prop}

\begin{itemize}
\item [1.] $ \gin {\ga_{K_X}(V)} = \gin V \, (\mbox{over } K_X). $
\item [2.] $ \gin { \ga_{K_X}(V)_{|x_s = 0}} = \gin V _{|x_s = 0} \, 
 (\mbox{over } K_X).$

\end{itemize}

\end{prop}

\begin{proof} We prove 2. The proof of 1. is analogous and easier. Besides
we will not need 1. We just state it for completeness.

\vskip 2mm
a) Let $S_1^{\si} = (x_1, \ldots, x_{s-1})$ and $G^{\si} = GL(S_1^{\si \vee})$.
Let $k \pil K$ be a homomorphism of fields and let
$G_K^{\si} = GL(S_1^{\si \vee} \te_k K)$. 
Due to the naturally split inclusion $S_1^\si \sus S_1$, there is a diagram

\[ \begin{CD}
\sas_K  \times G_K^\si @>>> \sas_K \\
@VVV @VVV \\
\saf_K \times G_K^\si @>>> \saf_K 
\end{CD} \]
where the upper action is given by $(a,g) \mapsto g^{-1}.a$.
The lower map gives a $K$-algebra homomorphism 
\[ \ga^{\si} : K[x_1, \ldots, x_s] \lpil K[G_K^\si] 
\te_K K[x_1, \ldots, x_s]. \]
The upper map gives a $K$-algebra homomorphism
\[ \ga^{\si}_{|x_s = 0} : K[x_1, \ldots, x_{s-1}] \lpil K[G_K^\si] \te_K
 K[x_1, \ldots, x_{s-1}]. \]
For a homogeneous subspace $W \sus K[x_1, \ldots, x_s]$ we now see that
\[ \ga_{|x_s = 0}^{\si}(W_{|x_s = 0}) = \ga^{\si}(W)_{|x_s = 0}. \]
The initial space of the former is by definition
$\gin {W_{|x_s = 0}}$.
By 1.3
applied to the latter initial space we then get
\begin{equation} \gin {W_{|x_s = 0}} = \ii { \ga^{\si}(W)}_{|x_s = 0}.
\label{li1} \end{equation}

\vskip 3mm

b) Now there is a diagram

\[ \begin{CD}
\saf \times G^\si \times G @>>> \saf \times G \\
@VVV @VVV \\
\saf \times G @>>> \saf. 
\end{CD} \]

The upper horizontal map is given by
$(a,h,g) \mapsto (h^{-1}.a,g)$. The lower horizontal map and the right
vertical map are the actions. Lastly, the left vertical map is given by
$(a,h,g) \mapsto (a,hg)$.
It induces a diagram

\begin{equation} \begin{CD}
\saf \times G^\si \times X @>>> \saf \times X \\
@VVV @VVV \\
\saf \times G @>>> \saf. 
\end{CD}  \label{li2} \end{equation}
Apply lemma 2.7.
Then $\ga_{K_X}(V)$ has initial space $\gin V$. Also applying 2.7
to the composition $\saf \times G^\si \times X \pil \saf \times G \pil \saf$
(from the diagram), gives that $\ga_{K_{G^\si \times X}} (V)$ has 
initial ideal $\gin V$ over $K_{G^\si \times X}$.


Now go back to part a) of this proof and put $K = K_X$ and 
$W = \ga_{K_X} (V)$. By the commutativity of the diagram (\ref{li2}) we see
that
\[ \ga^{\si} (W) = \ga_{K_{G^\si \times X}} (V). \]
Thus 
\[ \ii {\ga^{\si} (W)} = \ii {\ga_{K_{G^\si \times X}} (V)}
= \gin V \, (\mbox{over } K_{G^\si \times X}). \]
Putting this together with (\ref{li1}) we get
\[ \gin { \ga_{K_X}(V)_{|x_s = 0}} = \gin V _ {|x_s = 0} \,
   (\mbox{over } K_X). \]
\end{proof}

\parg Now, there is of course also a natural action
\[ G \times \saf \lpil \saf \]
given by 
\[ (g,a) \mapsto g.a. \]
The morphism 
\[ \rho : G \times \saf \lpil \saf \times G \]
given by
\[ (g,a) \mapsto (g.a, g)\]
is an isomorphism and its inverse $\rho^{-1}$ is given by
\[ (b,g) \mapsto (g, g^{-1}.b),\]
The morphism $\rho$ induces a $k[G]$-algebra isomorhism
\[ \Ga : k[x_1, \ldots, x_s] \te_k k[G] \lpil k[G] \te_k k[x_1, \ldots, x_s].\]
Note that $\Ga^{-1}$ is the $k[G]$-algebra isomorphism induced by $\rho^{-1}$.
For any $k[G]$-algebra $R$ we get an $R$-algebra isomorphism
\[ \Ga_R : R[x_1, \ldots, x_s] \lpil R[x_1, \ldots, x_s]. \]

The homogeneous subspace $V \sus S$ induces an $R$-submodule
\[ V_R = V \te_k R \sus R[x_1, \ldots, x_s] \]
and so we get a free $R$-module
\[ \Ga_R^{-1}(V_R) \sus R[x_1, \ldots, x_s] \]
which is in fact just $\ga_R(V)$.

\parg For a morphism $\phi : X \pil G$ with $\phi(X) \cap U \neq 0$ we now
see that 
\[ \ii { \Ga_{K_X}^{-1}(V_{K_X})} = \ii {\ga_{K_X} (V)} 
   = \gin V \, (\mbox{over } K). \]
Now consider $G = GL(S_1^{\vee})$ to be an open subset of ${\bf A}^{s^2}$
with coordinate functions $u_{ij}$ for $i,j = 1, \ldots, r$. 
Let the $u_{ij}$ take general values of $k$ for $i < r$ and 
let $u_{rj} = t_j$. Let $D$ be the determinant of the matrix thus obtained
and let $T = k[t_1, \ldots, t_s]_D$. The situation to which we will apply
the above is to the situation where $X = $ Spec $T$. For the rest
of the paper let $K = K_X = k(t_1, \ldots, t_s)$ the field of rational 
functions in the $t_i$.

Finally, if we let $h = \sum_{i=1}^s t_i x_i$, note the following will be 
used repeatedly in sections 4 and 5 : $\Ga_K (x_s) = h$.

\section{ Derivatives of forms}

Given a form $p$ in $S_K = K[x_1, \ldots, x_s]$.
One may then differentiate it with respect to the $t_i$ and obtain partial
derivatives $\pd^{|I|}p / \pd t^I$ where $I = (i_1, \ldots, i_r)$.
More generally for a homogeneous form 
$s(\bt) = \sum \alpha_I t^I$ of degree $d$ we get the directional derivative
$\pd^d p / \pd s(\bt) = \sum \alpha_I \pd^d p/ \pd t^I$.

For a form $f$ in $S_K$ let $\ba f$ be its image in 
$S_K / (h)$.
Now consider a specific form $p$. Let $l(\bt) = \sum \alpha_i t_i$ 
be such that 
$\ba {l(\bx)} = \sum \alpha_i \ba {x_i}$ is not a factor of $\ba p$.

\begin{lemma} Suppose $\pd^k p/ \pd l^k = \alpha_k p$ for $k \geq 0$.
If $f$ is a form such that $\ba {\pd^k f / \pd l^k}$ has $\ba p$
as a factor for all $k \geq 0$, then $f$ has $p$ as a factor.
\end{lemma}

\begin{proof} We have 
\begin{equation} f = u_1 p + h a_1 \label{a1} \end{equation}
for some $u_1$ and $a_1$.
 Differentiating this gives
\[ \pd f/ \pd l = \pd u_1 / \pd l \cdot p + u_1 \pd p / \pd l + l(\bx) a_1
+ h \pd a_1 / \pd l. \]
Thus $\ba p$ divides $\ba a_1$. So $a_1 = v_1 p + h a_2$ for some $v_1$ and 
$a_2$. Inserting this in (\ref{a1}) gives
\[ f = u_2 p + h^2 a_2  \]
where $u_2 = u_1 + h v_1. $
Now differentiate twice with respect to $l$. We may conclude that
\[ a_2 = v_2 p + h a_3 \]
for som $v_2$ and $a_3$.
Continuing we get in the end that $f = up$.
\end{proof}

The following result is proposition 10 in \cite{Co} and is due to Green. 
It is assumed there that the field $k = {\bf C}$ but the proof is readily
seen to work for any field of characteristic zero.
Given a form $\ba p$ in $S_K / (h)$ it
gives a criterion for it to lift to a form in $S_K$ which
is essentially a form in $S$.

\begin{prop} Let $p \in S_K$ be a form such that
\[ x_i \ba{ \pd p / \pd t_j} \equiv x_j \ba{ \pd p / \pd t_i} \pmod {\ba p} \]
for all $i$ and $j$.
Then $p = \alpha p_0 + h R$ where $p_0 \in S$ and
$\alpha \in K$.
\end{prop}

Consider now a form $f \in S \sus S_K$.
It gives a hypersurface in $\bf P^{s-1}$. The following says that if all
hyperplane sections of this hypersurface are reducible with a component of
a given degree then the same is true for the hypersurface defined by $f$.

\begin{kor} Suppose $\ba f = \ba u \cdot \ba p$ in $S_K / (h)$,
where $\ba u$ and $\ba p$ do not have a common factor. Then $\ba p$ lifts 
to a form $\alpha p_0$ where $p_0 \in S$.
Furthermore $p_0$ is a factor of $f$.
\end{kor}

\begin{proof} Let $u$ and $p$ in $S_K$ be liftings of 
$\ba u$ and $\ba p$. We get 
\[ f = up + hR. \]
Differentiating with respect to $\pd / \pd t_i$ gives
\[ 0 = \pd u / \pd t_i \cdot p + u \pd p / \pd t_i 
  + x_i R + h \pd R / \pd t_i. \]
Thus we get 
\[ \ba u ( x_j {\ba {\pd p/ \pd t_i}} - x_i \ba {\pd p / \pd t_j})
   \equiv 0  \pmod {\ba p}. \]
Then by proposition 3.2 we conclude that $\ba p$ has a lifting
$\alpha p_0$ where $p_0 \in S$. 
By lemma 3.1 we conclude that $p_0$ is a factor of $f$ since the
$\pd^k f / \pd l^k = 0$ for $k \geq 1$.
\end{proof}

   Now suppose $V \sus S_{n+m}$ is a subspace so we get a subspace
$V_K = V \te_k K \sus S_{K,n+m}$
and $V_{K|h=0} \sus S_K / (h)$.

\begin{prop} Suppose the forms of $V_{K|h=0}$ have a common factor $\ba p$
where $\ba p$ is a common factor of maximal degree $m$. Then $V$ has 
a common factor $p_0$ of degree $m$ such that $\ba p = \alpha \ba p_0$
for some $\alpha \in K$.
\end{prop}

\begin{proof}
We may choose an $f_0 \in V$ such that
\[ \ba f_0 = \ba u_0 \ba p \]
where $\ba u_0$ and $\ba p$ are relatively prime.
This is seen as follows. Let $\ba p = {\ba a_1}^{e_1} \cdots {\ba a_r}^{e_r}$
be a factorization where the $\ba a_i$ are distinct irreducible factors.
It is easily seen that the set of $f$ in $V$ where $\ba f$ has 
${\ba a_i}^{e_i + 1}$ as a factor, is a linear subspace $V_i$ of $V$.
On the other hand if $f$ varies all over $V$ the restrictions $\ba f$ 
generate $V_{K|h=0}$.
Thus we cannot have $V_i = V$ for any $i$.
But since char $k= 0$ the field $k$ is infinite, so there must be an $f_0$ in 
$V - \cup V_i$.


By corollary 3.3, $\ba p$ lifts to $\alpha p_0$ where $p_0 \in S$.
Choose now any $f$ in $V \sus V_K$. Then 
\[ \ba f = \ba u \cdot \ba {\alpha p_0}. \]
By lemma 3.1 we may conclude that $p_0$ is a factor of $f$ and thus a common
factor of $V$.
\end{proof}

\section{ The case when $s = r$}

Now we are ready for the specific work in proving the Main Theorem.
Consider $S = k[x_1, \ldots, x_r]$. Let $W = (x_1, \ldots, x_r) = S_1$ which
is a linear space. Use the notation $W^n = S_n$. (This will make our 
statements more unified in form.) 
Let the monomial order be the revlex order.
In this section we prove the following (which
is the case $s=r$ of the Main Theorem.)

\begin{theorem} Let $V \sus S_{n+m}$ be a linear space such that 
\[ \gin V = W^n x_1^m \sus S_{n+m}. \]
Then there exists a polynomial $p \in S_{m}$ such that $V = W^n p.$
\end{theorem}

We assume $V$ to be in general coordinates so 2.7
applies.

\begin{lemma} There is a form $p$ in $S_{K,m}$ such that
\[ h^n p \in V_K. \]
\end{lemma}

\begin{proof}
From 2.11 we have $\ii {\Ga_K^{-1}(V_K)} = \gin V$ over $K$.
Thus there exists a $q_0$ in $\Ga_K^{-1}(V_K)$ such that
\[ q_0 = x_r^n x_1^m + \mbox{terms with smaller monomials}.\]
By the property of the revlex order, $x_r^n$ will divide
all terms of $q_0$ so there exists a $p_0 \in S_{K,m}$ such that
\[ q_0 = x_r^n p_0. \]
Let $p = \Ga_K (p_0)$. Then we get
\[ h^n p = \Ga_K(x_r)^n \Ga_K(p_0) = \Ga_K(q_0) \in V_K. \]
\end{proof}

From $V_K \sus S_K$ we obtain the subspace
\[ V_{K|h=0} \sus S_K / (h). \]
Let $\ba p$ be the image of $p$ in $V_{K|h=0}$.

\begin{lemma} The elements in $V_{K|h=0}$ have $\ba p$ as a common factor.
Furthermore it is a common factor of maximal degree. 
\end{lemma}

\begin{proof} We first find the dimension of the space $V_{K|h=0}$.
The map $\Ga_K$ gives an isomorphism
\[ \ba {\Ga_K} : K[x_1, \ldots, x_r]/ (x_r) \lpil K[x_1, \ldots, x_r]/ (h). \]
Thus $\ba {\Ga_K}^{-1}(V_{K|h=0}) = \Ga_K^{-1}(V_K)_{|x_r = 0}$.
Since $\Ga_K^{-1}(V_K)$ has initial space
\[ (x_1, \ldots, x_r)^n \cdot x_1^m, \]
we get by 1.3 that $\Ga_K^{-1}(V_K)_{|x_r = 0}$ has initial space
\[ (x_1, \ldots, x_{r-1})^n \cdot x_1^m. \]
Hence the dimension of $V_{K|h=0}$ is equal to the dimension of this space.

Now differentiate the equation 
\[ h^n p \in V_K \]
with respect to $ \partial^{|I|}/ \partial t^I$ where 
$I = (i_1, \ldots, i_{r-1})$
and $|I| = n$.
The derivative will also be in $V_K$. This is essentially the fact that
when a vector varies in a vector space the derivatives will also 
be in that vector space.
We thus get
\[ \bx^I p + h R_I \in V_K \]
for some $R_I$. Thus 
\begin{equation} \bx^I \ba p \in V_{K|h=0}. \label{bap} \end{equation}
But when $I$ varies, all these forms are linearly independent
since $h$ does not divide any linear combination of the $\bx^I$.
By our statement about the dimension of $V_{K|h=0}$, the forms
(\ref{bap}) must generate $V_{K|h=0}$, thus proving the lemma.
\end{proof}

By corollary 3.4 we may now conclude that $V$ has a maximal common factor
$p_0$ of degree $m$. Thus proving 4.1.

\section {The case when $s > r$}

Now we assume $S = k[x_1, \ldots, x_s]$. As before 
$W = (x_1, \ldots, x_r) \sus S_1$, a linear subspace and assume $s > r$.
The monomial order is revlex.
In this section we prove the following by induction on $s$.

\begin{theorem} Let $V \sus S_{n+m}$ be a linear space such that 
\[ \gin V = W^n x_1^m \sus S_{n+m}. \]
Then there exists a polynomial $p \in S_{m}$ and a linear subspace
$W_n \sus S_n$ such that $V = W_n p.$
\end{theorem}

Assume $V$ to be in general coordinates.
Let $g : S_1 \pil S_1$ be a general coordinate change. Since
$\ii {g^{-1}.V} = (x_1, \ldots, x_r)^n \cdot x_1^m$, by 1.3 it follows that 
$ \ii {g^{-1}.V_{|x_s = 0}} = (x_1, \ldots, x_r)^n \cdot x_1^m$ also.
By induction $g^{-1}.V_{|x_s = 0}$ has a common factor. By translating back, 
$V_{|g.x_s = 0}$
also has a common factor (depending on $g$). 
The following expresses this in the algebraic language we use.

\begin{lemma} There is a form $p$ in $S_{K,m}$ such that $\ba p$ in 
$S_{K|h=0}$ is a common factor of $V_{K|h=0}$. Furthermore it is
a common factor of maximal degree.
\end{lemma}

\begin{proof} By 2.9.2 the generic initial ideal of 
$\Ga_K^{-1}(V_K)_{|x_s = 0} = \ga_K(V)_{|x_s = 0}$ is 
$\gin V _{|x_s = 0}$ (over $K$).
The latter is, by 1.3, seen to be
\[ (x_1, \ldots, x_r)^n \cdot x_1^m. \]
By induction there is a form $\overline{p_1}$ in $S_{K,m|x_s = 0}$ which is
a common factor of $\Ga_K^{-1}(V_K)_{|x_s = 0}$. 
Now $x_1^m$ is a common factor of $\ii {\Ga_K^{-1}(V_K)_{|x_s = 0}}$
of maximal degree. Then $\ba p_1$ must also have maximal degree
as a common factor of $ \Ga_K^{-1}(V_K)_{|x_s = 0}$.
Lift this to a form $p_1$ in $S_{K,m}$. Then $p= \Ga_K(p_1)$ is the 
required form.
\end{proof}

   By corollary 3.4 we may now conclude that $V$ has a maximal common factor 
$p_0$ of degree $m$. Thus proving 5.1.

\section {An example}

   Consider the complete intersection of three quadratic forms in $\ptre$.
Let $I \sus k[x_1,x_2,x_3,x_4]$ be its homogeneous ideal. We have the following
facts.

\begin{itemize}
\item [1.] $I$ and $\gin I$ have the same Hilbert functions.

\item [2.] $\gin I$ is Borel-fixed. ( See proposition 15.20 in \cite{Ei}.)

\item [3.] Since $I$ is saturated, by proposition 2.21 in \cite{Gr} we have 
      $\gin I : x_4 = \gin I$. This is really just the fact that 
   $\ii {I : x_4} = \ii I : x_4$ for the revlex order
   (proposition 15.12 b. in \cite{Ei}),   
   and that if $I$ is in general coordinates and saturated then 
   $I : x_4 = I$.

\end{itemize}

These three facts imply that there are two possible candidates for $\gin I$~:


\begin{eqnarray*} 
J^{(1)} & = & (x_1^2, x_1x_2, x_2^2, x_1x_3^2, x_2x_3^2, x_3^4), \\
 J^{(2)} & = & (x_1^2, x_1x_2, x_1x_3, x_2^3, x_2^2x_3, 
x_2x_3^2, x_3^4).
\end{eqnarray*}

However, by the theorem above if $\gin I = J^{(2)}$ then the quadratic forms
in $I_2 \sus S_2$ would have to have a common factor. Impossible.
Thus $\gin I = J^{(1)}$.  On the other hand, if $I$ is an ideal with $\gin I = J^{(2)}$
then since the quadratic forms in $I_2$ would have a common factor it must 
be the ideal of seven points in a plane pluss one extra point not in the
plane.

Note also the following. Let $>_1$ be the ordering of the monomials which is
lexicographic in the three first variables, and then refined with the 
reverse lexicographic order with respect to the last variable. I.e.
\[ x_1^{a_1}x_2^{a_2}x_3^{a_3}x_4^{a_4} > x_1^{b_1}x_2^{b_2}x_3^{b_3}x_4^{b_4}
\] 
if $a_4 < b_4$, or $a_4 = b_4$ and
\[  x_1^{a_1}x_2^{a_2}x_3^{a_3} > x_1^{b_1}x_2^{b_2}x_3^{b_3}
\]
for the lexicographic order.
Then if the three forms are general it is easily seen that $\gin I = J^{(2)}$. In 
fact it is not difficult to argue that one will always have $\gin I = J^{(2)}$
if you have a complete intersection of three forms and this order.
Thus both $J^{(1)}$ and $J^{(2)}$ are in fact specialisations of $I$.


Furthermore it is not difficult to give an example of a complete intersection
of three forms such that in$(I) = J^{(2)}$ for the reverse lexicographic order.
Thus the fact that one can read some interesting algebraic or geometric
information from the initial ideal depends on the fact that you are looking
at the {\it generic initial ideal}.

To sum up, $J^{(2)}$ is a specialisation of the ideal $I$ of a complete 
intersection of three quadratic forms in general coordinates through
the order $>_1$ given above. It is also the specialisation of an ideal $I$
of a complete intersection of three quadratic forms through the
revlex order, but it is {\em never} a specialisation of the ideal $I$ 
of a complete intersection of three quadratic forms through the revlex
order when the forms are in general coordinates.

\end{document}